\documentstyle[12pt]{article}

\begin{document}



\begin{center}

{\bf FOURIER TRANSFORMS IN EXPONENTIAL  REARRANGEMENT INVARIANT SPACES} \\

\vspace{3mm}

$ {\bf E.Ostrovsky^a, \ \ L.Sirota^b } $ \\

\vspace{2mm}

$ ^a $ Corresponding Author. {\it Department of Mathematics, Ben Gurion University of Negev,
84105, Beer - Sheva, Ben Gurion street, 2, Israel.}\\
\end{center}
E - mail: \ galaostr@cs.bgu.ac.il\\
\begin{center}
$ ^b $ {\it Department of Mathematics. Bar Ilan University, 59200, Ramat Gan, Israel.}\\
\end{center}
E - mail: \ sirota@zahav.net.il\\

\vspace{2mm}

{\bf Abstract}.  In this article we investigate the Fourier series and transforms for 
the functions defined on the $ [0, 2 \pi]^ d $ or $ R^d $ and  belonging to  the exponential 
Orlicz and some other rearrangement invariant (r.i.) spaces.\\



{\it Key words: }
 Orlicz and rearrangement invariant (r.i.) spaces, Fourier Integrals and Series, 
Operators, Moment Inequalities, fundamental Function, equivalent Norms, slowly varying 
Functions, Wavelets, Haar's Series.



{\it  Mathematics Subject Classification 2000.} Primary 42Bxx, 4202; Secondary 28A78, 42B08. \\
\section{Introduction. Notations. Problem Statement.} 
 For the (complex) measurable function $ f = f(x) $ defined on the $ X = \{x\} = T^d = [0, \ 2 \pi]^d, 
\ d = 1,2, \ldots $ or $ X = R^d $ we denote correspondently the Fourier coefficients and 
transform
$$
c(n)=\int_{T^d} \exp( i(n,x) ) \ f(x) \ dx, \ \ F[f](t) = \int_{R^d} \exp(i(t,x)) \ f(x) \ dx,
$$
 where as usually
$$
F[f](t) \stackrel{def}{=}  \lim_{M \to \infty}\int_{|x| \le M} \exp(i(t,x)) \ f(x) \ dx,
$$
$$
x = (x_1,x_2,\ldots,x_d), \  n = (n_1,n_2, \ldots, n_d),  \ t = (t_1,t_2, \ldots, t_d), \ 
dx = \prod_{j=1}^d dx_j
$$
if $ x \in R^d, $ and $ dx = (2 \pi)^{-d} \ \prod_{j=1}^d dx_j $ in the case $ X = T^d; $
$$
n_i = 0, \pm 1, \pm 2, \ldots, \ (t,x) = \sum_{j=1}^d t_j x_j, \ |n| = \max_{j} |n_j|, \ 
|t| = \max_j |t_j|, 
$$
$$
s_M[f](x) = (2 \pi)^{-d}  \sum_{n: |n| \le M} c(n) \exp(-i (n,x)), \ x \in T^d; 
$$
$$
 S_M[f](x) = (2 \pi)^{-d} \int_{t: |t| \le M} \exp(-i(t,x)) \ F[f](t) \ dt, \ x \in R^d. 
$$
 We will assume in this paper that the function $ f(\cdot) $ belongs to some Orlicz space 
$ L(N) = L(N;X) $   with so - called {\it exponential 
$ N \ - $ function } $ N = N(u) $ or some rearrangement invariant  space $ E $ with fundamental 
function $ \phi_E(\delta), \ \delta > 0, $ and will investigate the properties of Fourier 
transform of $ f, $ for example, the boundness of operators $ S_M[\cdot], \ s_M[\cdot] $ and 
the convergence and divergence in {\it some} Orlicz norms
 $ s_M[f](\cdot) \to f, \ S_M[f](\cdot) \to f  $ as $ M \to \infty. $  \par
  Note than the case if the function $ N(\cdot) $ satisfies the $ \Delta_2 $ condition is 
known; see, for example, [13], [11, p 193 - 197]. Our results are also some 
generalization of [4]. The case of Hilbert transform and generalized functions 
$ \ \{ f \} \ $ see in [2.] \par
 Now we will describe the using Orlicz and other  spaces. Let $ (X,A,\mu) $ be some 
measurable space with non - trivial measure $ \mu. $ For the measurable complex function
$ f(x) = f_1(x) + i f_2(x), x \in X, f_{1,2}: X \to R^1 $ the symbol
$ |f|_p = |f|_p(X,\mu) $ will denote the usually $ L_p $ norm:
$$
|f|_p = ||f||L_p(X, \mu) = \left[ \int_X |f(x)|^p \ \mu(dx) \right]^{1/p}, \ p \ge 1.
$$
In the case $ X = R^d $ we introduce a new measure $ \nu(\cdot) $ (non - finite, in 
general case): for all Borel set $ A \subset R^d $
$$
\nu(A) = \int_A \prod_{i=1}^d x_i^{-2} \ dx = \int_A \prod_{i=1}^d x_i^{-2} \ \cdot 
\prod_{i=1}^d dx_i,
$$
and will denote  $ |f|_p(\nu) = $
$$
\left| \left| \left(\prod_{i=1}^d x_i \right) \cdot f \right| \right|L_p(X,\nu) = 
\left[ \int_X \left|\prod_{j=1}^d x_j \right|^p \ \cdot |f(x)|^p \ \nu(dx)  \right]^{1/p} =
$$
$$
\left[\int_X \left| \prod_{j=1}^d x_j \right|^{p-2} \cdot |f(x)|^p  \ dx    \right]^{1/p}.
$$
 For arbitrary multiply sequence (complex, in general case)
$ c(n) = c(n_1,n_2,\ldots, n_d), \\ n_i = 0, \pm 1, \pm 2, \ldots,  n \in Z^d $ we denote 
as usually 
$$
|c|_p  =  \left[ \sum_{n} |c(n) |^p \right]^{1/p}, \ p \ge 1;
$$
and introduce the discrete analog of $ |f|_p(\nu) $ norm:
$$
|c|_p(\nu) = \left[ \sum_n \left|c(n) \right|^p \ \cdot  \left(
\left| \prod_{j=1}^d n_j \right|^{p-2} +1 \right) \right]^{1/p}, \ p \ge 2.
$$
 Let $ N = N(u) $ be some 
$ N \ - $ Orlicz's function, i.e. downward convex, even, continuous 
differentiable for all sufficiently greatest  values $ u, \ u \ge u_0, $
strongly  increasing in the right - side axis,  
and such that $ N(u) = 0 \ \Leftrightarrow  u = 0; \ u \to \infty \ \ \Rightarrow 
dN(u)/du \ \to \ \infty. $ We say that  $ N(\cdot) $ is an Exponential Orlicz Function, 
briefly: $ N(\cdot) \in \ EOF, $ if $ N(u) $ has a view: for some continuous 
differentiable strongly increasing downward convex in the domain $ [2,\infty] $ 
function $ W = W(u) $ such that $ u \to \infty \ \Rightarrow W^/(u) \to \infty $

$$
N(u) = N(W,u) = \exp(W(\log |u|)), \ |u| \ge e^2.
$$
 For the values $ u \in [-e^2,e^2] $ we define $ N(W,u) $ arbitrary  but so that the function 
$ N(W,u) $ is even continuous convex strictly increasing in the right side axis
and such that $ N(u) = 0 \ 
\Leftrightarrow u = 0. $ The correspondent Orlicz space on $ T^d, \ R^d $ with usually 
Lebesque measure  with $ N \ - $ function $ N(W,u) $ 
we will denote $ L(N) = EOS(W); \ EOS = \cup_{W} \{EOS(W)\} $ (Exponential Orlicz's Space).\par
 For example, let $ m = const > 0, \ r = const \in R^1, $ 
$$
N_{m,r}(u) = \exp \left[|u|^m \ \left(\log^{-mr}(C_1(r) + |u| ) \right) \right] -1, 
$$
$ C_1(r) = e, \ r \le 0; \ C_1(r) = \exp(r), \ r > 0. $ Then $ N_{m,r}(\cdot) \in EOS. $ 
In the case $ r = 0 $ we will write $ N_m = N_{m,0}. $ \par 
 Recall here that the Orlicz's norm on the arbitrary measurable space $ (X,A,\mu) $
$ ||f||L(N) = L(N,X,\mu) $ may be calculated by the formula (see, for example,  
[11], p. 73; [6], p. 66) 
$$
||f||L(N) = \inf_{v > 0}  \left \{ v^{-1} \left(1 + \int_X N( v|f(x)|) \ \mu(dx) \right) \right \}.
$$
 Recall also that the notation $ N_1(\cdot) << N_2(\cdot) $ for two Orlicz functions 
$ N_1, N_2 $ denotes:
$$
\forall \lambda > 0 \ \Rightarrow \lim_{u \to \infty} N_1(\lambda u)/N_2(u) = 0.
$$
 We will denote for arbitrary Orlicz $ L(N) $ (and other r.i.) spaces by $ L^0(N) $ the 
closure of all bounded functions with bounded support. \par 
 Let $ \alpha $ be arbitrary number, $ \alpha \ge 1, $ and $ N(\cdot) \in EOS(W) $ for 
some $ W = W(\cdot). $ We denote for such a function $ N = N(W,u) $
by $ N^{(\alpha)}(u) $ a new $ N \ - $ Orlicz's function such that
$$
N^{(\alpha)}(u) = C_1 \ |u|^{\alpha}, \ \ |u| \in [0, C_2];
$$
$$
N^{(\alpha)}(u) = C_3 + C_4 |u|, \ \ |u| \in (C_2, C_5];
$$
$$
N^{(\alpha) }(u) = N(u), \ \ |u| > C_5, \ \ 0 < C_2 < C_5 < \infty,\eqno(1.1) 
$$
$$
C_{1,2,3,4,5} = C_{1,2,3,4,5}(\alpha,N(\cdot)).
$$
 In the case $ \alpha = m(j+1), \ m > 0, \ j = 0,1,2, \ldots $ the function $ N^{(\alpha)}_m(u) $
is equivalent to the following Trudinger's function:
$$
N_m^{(\alpha)}(u) \sim N_{[m]}^{(\alpha)}(u) = \exp \left(|u|^m \right) -\sum_{l=0}^{j}
u^{ml}/l!.
$$
 This method is described in [15], p. 42 - 47. Those Orlicz spaces are applied to the theory
of non - linear partial differential equations [15].\par
 We can define formally the spaces $ L(N^{(\alpha)}_m) $ at $ m = + \infty $ as a projective 
limit at $ m \to \infty $ the spaces $ L(N^{(\alpha)}_m), $ but it is evident that
$$
L \left(N^{(\alpha)}_{\infty} \right) \sim  L_{\alpha} + L_{\infty},
$$
where the space $ L_{\infty} $ consists on all the a.e. bounded functions with norm 
$$ 
|f|_{\infty} = vraisup_{x \in X} |f(x)|.
$$
 Of course, in the case $ X = T^d $ 
$$
L_{\alpha} + L_{\infty} \sim L_{\infty}.
$$
 Here and further we will denote by $  C_k = C_k(\cdot), k = 1,2,\ldots $ some positive finite   
essentially and by $ C, C_0 $ non - essentially "constructive" constants. 
 It is very simple to prove the existence of constants $ C_{1,2,3,4,5} = C_{1,2,3,4,5}(\alpha, 
N(\cdot)) $ such that $ N^{(\alpha)} $ is some new exponential $ N \ $
Orlicz's function.  By the symbols $ K_j $ we will denote the "classical" absolute constants.\par  
  Now we will introduce some {\it new } Banach spaces. Let $ \alpha = const \ge 1 $ and
$ \psi = \psi(p), \ p \ge \alpha $
be some continuous positive: $ \psi(\alpha) > 0 $  finite strictly increasing function such 
that the function $ p \to p \log \psi(p) $ is downward convex and 
$$
\lim_{p \to \infty} \psi(p) = \infty.
$$
The set of all those functions we will denote $ \Psi; \ \Psi = \{\psi\}. $ A particular case:

$$
\psi(p) = \psi(W;p) = \exp(W^*(p)/p), 
$$
where
$$
W^*(p) = \sup_{z \ge \alpha} (pz - W(z))  
$$
is so - called Young - Fenchel, or Legendre transform of $ W(\cdot). $ It follows from 
theorem of Fenchel - Moraux that in this case 
$$
W(p) = \left[ p \ \log \psi(W;p) \right]^*, \ \ p \ge p_0 = const \ge 2,
$$
and consequently  for all $ \psi(\cdot) \in \Psi $ we introduce the correspondent $ N \ - $ 
function by equality:
$$
N([\psi]) = N([\psi],u) =  
\exp \left \{\left[p \log \psi(p) \right]^*(\log u) \right \}, \ u \ge e^2.
$$
Since $ \ \forall \ \psi(\cdot) \in \Psi, \ d = 0,1,\ldots \ \Rightarrow p^d \cdot \psi(p) \in \Psi, $
we can denote
$$
\psi_d(p) = p^d \cdot \psi(p), \ \  N_d([\psi]) = N_d([\psi],u) = N([\psi_d], u). 
$$
For instance, if $ N(u)= \exp(|u|^m), \ u \ge 2, $ where $ \ m = const > 0, $ then
$$
N_d([\psi], u) \sim \exp \left(|u|^{m/(dm+1)} \right), \ u \ge 2.
$$
{\bf Definition 1. } We introduce for arbitrary such a function $ \psi(\cdot) \in \Psi $ 
the so - called $ G(\alpha;\psi) \ $ and $ G(\alpha; \psi,\nu) $
 norms  and correspondent Banach spaces $ G(\alpha; \psi), \ G(\alpha, \psi, \nu) $ as a 
set of all measurable complex functions  with finite norms:

$$
||f||G(\alpha; \psi) = \sup_{p \ge \alpha} (|f|_p/\psi(p)),
$$
and  analogously,

$$
||f||G(\alpha; \psi,\nu) = \sup_{p \ge \alpha} (|f|_p(\nu)/\psi(p)).
$$
 For instance $ \psi(p) $ may be $ \psi(p) = \psi_m(p) = p^{1/m}, \ m = const > 0; $ in 
this case we will write $ G(\alpha, \psi_m) = G(\alpha,m) $ and
$$
||f||G(\alpha,m) = \sup_{p \ge \alpha}  \left(|f|_p \ p^{-1/m} \right).
$$
 Also formally  we define
$$
||f||G(\alpha, m) = |f|_{\alpha} + |f|_{\infty}.
$$
{\it Remark 1.} It follows from Iensen inequality that in the case $ X = T^d $ all the spaces 
$ G(\alpha_1; \psi), \ G(\alpha_2,\psi), 1 \le \alpha_1 < \alpha_2 < \infty $ are isomorphic: 
$$
||f||G(\alpha; \psi) \le ||f||G(1; \psi) \le \max(1, \psi(\alpha)) \ ||f||G(\alpha; \psi)
$$
 It is false in the case $ X = R^d.$ \par
{\it Remark 2. } $ G(\alpha; \psi) $ is a rearrangement invariant (r.i.) space. $ G(\alpha,m) $
has a  fundamental function $ \phi(\delta; G(\alpha,m)), \ \delta > 0, $ where for any 
rearrangement invariant space $ G $
$$
\phi(\delta;G) \stackrel{def}{=} ||I_A(\cdot)||G(\cdot), \ \  mes(A) = \delta \in (0,\infty),
$$
$ mes(A) $ denotes usually Lebesque measure of Borel set $ A. $  We have:
$$
 \phi(\delta; G(\alpha,m)) = (e \ m \ |\log \delta|)^{-1/m}, \ \ \delta > \exp(-\alpha/m),
$$
$$
\phi(\delta; G(\alpha,m)) = \alpha^{-1/m} \ \delta^{1/\alpha}, \ \delta \in (0, \exp(-\alpha/m)).
$$
Let us consider also another space 
$ G(a,b,\alpha, \beta), \ 1 \le a < b < \infty; \ \alpha,\beta \ge 0. $ Here
$ X = R^d $ and we denote $ h = \min((a+b)/2; 2a). $  We introduce the function 
$ \zeta: (a,b) \to R^1_+: $ 
$$
 \zeta(p)= \zeta(a,b,\alpha,\beta; p) = (p-a)^{\alpha}, \ p \in (a, h); 
$$

$$
 \zeta(p) = (b-p)^{\beta}, \ p \in [h, b).
$$
 {\bf Definition 2.} The space $ G(a,b,\alpha,\beta) $ consists
on all the complex measurable functions  with finite norm:
$$
||f||G(a,b,\alpha,\beta) = \sup_{p \in (a,b)} \left[|f|_p \ \cdot \zeta(a,b,\alpha, \beta; p) 
\right].
$$
 The space $ G(a,b,\alpha,\beta) $ is also a rearrangement invariant space. After some 
calculations we receive the correspondence fundamental function
$ \phi(\delta; G(a,b,\alpha,\beta)). $ Namely, let us denote
 $ \delta_1 =  \exp(\alpha h^2/(h-a) ) $ and define for $ \delta \ge \delta_1 $
$$
p_1 = p_1(\delta) = \frac{\log \delta}{ 2 \alpha} - \left[\frac{\log^2 \delta}{4 \alpha^2} -        
\frac{a \log \delta}{\alpha} \right]^{1/2},
$$

$$
\phi_1(\delta) = \delta^{1/p_1} \ (p_1 - a)^{\alpha}; 
$$
and for $ \delta \in (0, \delta_1) \ \Rightarrow $
$$
\phi_1(\delta) = \delta^{1/h} (h-a)^{\alpha}. 
$$
 Further, set $ \delta_2 = \exp(-h^2 \beta/(b-h)), $ and for $ \delta \in (0, \delta_2) $ 
$ p_2 = p_2(\delta) = $
$$
- \frac{|\log \delta|}{2 \beta} + 
\left[ \frac{\log^2 \delta}{4 \beta^2} + b \frac{|\log \delta|}{\beta} \right]^{1/2},
$$

$$
\phi_2(\delta) = \delta^{1/p_2(\delta)} (b-p_2(\delta))^{\beta},
$$
and for $ \delta \ge \delta_2  $ we define
$$
\phi_2(\delta) = \delta^{1/h} (b-h)^{\beta}.
$$
We can write:
$$
\phi(\delta; G(a,b,\alpha,\beta)) = \max \left[ \phi_1(\delta), \phi_2(\delta) \right].
$$
 Note than at $ \delta \to 0+ $
$$
\phi(\delta; G(a,b,\alpha,\beta)) \sim \max \left \{ \delta^{1/h}(h-a)^{\alpha}, 
(e \beta b^2)^{\beta} \ \delta^{1/\beta}\ |\log \delta|^{-\beta} \right\}
$$
and at $ \delta \to \infty  $
$$
\phi(\delta; G(a,b,\alpha,\beta)) \sim \max \left\{\delta^{1/h} (b-h)^{\beta}, 
(2a^2 \alpha e^{-2})^{\alpha} \delta^{1/a} \ (\log \delta) ^{-\alpha} \right\}.
$$
 For example, let us consider the function $ f(x) = f(a,b; x), \ x \in R^1 \to R: $
$ f(x) = 0, \ x \le 0; $
$$
f(x) = x^{-1/b}, \ x \in (0,1); \ \ f(x) = x^{-1/a}, \ x \in [1,\infty); 
$$
then $ f(a,b,\cdot) \in G(a,b,1,1) $ and 
$$ 
\forall \ \Delta \in (0, 1/2] \ \Rightarrow f \notin G(a,b,1-\Delta,1) \cup G(a,b, 1, 1-\Delta). 
$$ 

 Analogously may be defined the "discrete" $ g(a,b,\alpha,\beta) $ spaces. Namely, let
$ c = c(n) = c(n_1,n_2,\ldots,n_d) $ be arbitrary multiply (complex) sequence. We say that 
$ c \in g(a,b,\alpha,\beta) $  if 

$$
||c||g(a,b, \alpha, \beta) \stackrel{def}{=} \sup_{p \in (a,b)}
\left[ |c|_p \ (p-a)^{\alpha} \ (b-p)^{\beta} \right].
$$
 It is evident that the non - trivial case of those spaces is only if $ \beta = 0; $ 
in this case we will write $ g(a,b, \alpha, 0) = g(a, \alpha) $ and 
$$
||c||g(\alpha) = \sup_{p > a} |c|_p \ (p-a)^{\alpha}.
$$
 We denote also for $ \psi(\cdot) \in \Psi: \ ||c||g(\psi, \nu) = $ 
$$
 \sup_{p \ge 2} \left[|c|_p(\nu)/\psi(p) \right], \ \ ||c||_m(\nu) = 
 \sup_{p \ge 2} \left[|c|_p(\nu) \ \cdot p^{-1/m} \right], \ \ m = const > 0.
$$

{\bf Our goal is investigating the boundness of Fourier operators and convergence (divergence)
 Fourier series and integrals in exponential Orlicz and $ G(a,b,\alpha, \beta),
\ \ g(\psi, \nu)  $ etc. spaces. } \\

 Note than our  Orlicz $ N \ - $ functions $ N \in EOS $ does not satisfy the so - called 
$ \Delta_2 $ condition. \\ 
\section{ Formulations of main results.}
{\bf Theorem 1.} {\it Let $ X = [0, 2 \pi]^d  $ and  $ \psi \in \Psi. $ Then the Fourier 
operators $ s_M[\cdot] $ are uniformly bounded in the space $ L(N[\psi]) $ into 
the other Orlicz's  space} $ L(N_d[\psi]): $

$$
\sup_{M \ge 1} ||s_M[f]||L(N_d[\psi]) \le C_6(d,\psi) \ ||f||L(N[\psi]). \eqno(2.1)
$$

{\bf Theorem 2.} {\it Let now $ X = R^d, \ \psi \in \Psi $ and $ \alpha = const > 1. $
The Fourier operators $ S_M[\cdot] $ are uniformly bounded in the space $ L(N^{(\alpha)}[\psi]) $
into the space } $ L(N^{(\alpha)}_d[\psi]): $
$$
\sup_{M \ge 1} ||S_M[f]||L(N^{(\alpha)}_d [\psi]) \le C_7(\alpha,d, \psi) \ 
||f||L(N^{(\alpha)}[\psi]). \eqno(2.2)
$$
 Since the function $ N[\psi] $ does not satisfies the $ \Delta_2 $ condition, the assertions 
(2.1) and (2.2) does not mean that in general case when $ f \in L(N_d^{\alpha}[\psi]) $  
$$
\lim_{M \to \infty} ||s_M[f] - f||L(N_d[\psi]) = 0, \eqno(2.3)
$$
$$
\lim_{M \to \infty} ||S_M[f] - f||L(N_d^{(\alpha)}[\psi]) = 0; \eqno(2.4) 
$$
see examples further. But it is evident that propositions (2.3) and (2.4) are true  if 
correspondently
$$ 
f \in L^0(N_d[\psi]), \ \ f \in L^0(N_d^{(\alpha)} [\psi]).
$$
Also it is obvious that if $ f \in L(N_d[\psi]), X = [0, 2 \pi]^d $ or, in the case $ X = R^d,
\ f \in L(N_d^{(\alpha)}[\psi]), $ then for all EOF $ \Phi(\cdot) $  such that $ \Phi << 
N_d[\psi] $ or $ \Phi << N_d^{(\alpha)}([\psi]) $  the following impications hold:
 $$ 
 \forall f \in L(N_d[\psi]) \ \Rightarrow \lim_{M \to \infty} ||s_M[f] - f||L(\Phi) = 0, \ X = 
[0,2\pi]^d; \eqno(2.5) 
 $$

$$
\forall f \in L(N_d[\psi]) \ \Rightarrow \lim_{M \to \infty} ||S_M[f]-f||L(\Phi) = 0, \ X = R^d.
\eqno(2.6)
$$
{\bf Theorem 3.} {\it Let $ \Phi(\cdot) $ be an EOF and let 
$ N(\cdot) =  L^{-1}(u) \in EOF, $ where $ L(y), \  y \ge \exp(2) $ is a positive slowly 
varying at $ u \to \infty $ strongly increasing continuous differentiable 
in the domain $ [\exp(2), \infty) $ function such that the function
$$
W(x) = W_L(x) = \log L^{-1}(\exp x), \ \ x \in [2, \infty)
$$
is again strong increasing to infty together with the derivative $ dW/dx. $ 
 In order to the impication (2.5) or,
correspondently, (2.6) holds, it is necessary and sufficient that $ \Phi << L(N_d[\psi]), $
or, correspondently } $ \Phi << L(N_d^{(\alpha)}[\psi]). $ \par
 For instance, the conditions of theorem 3 are satisfied for the functions $ N = N_{m,r}(u). $ \par
{\bf Theorem 4. } {\it Let $ f(\cdot) \in G(1,b, \alpha, 0), \alpha > 0. $ 
Then $ F[f] \in L(N^{(2)}_{1/\alpha}) $ and }
$$
\sup_{M \ge 1} ||S_M[f]||L \left(N^{(2)}_{1/\alpha} \right) \le C_8(\alpha, N) \ 
||f||G(1,b, \alpha, 0).  \eqno(2.7)
$$
{\bf Theorem 5. } {\it Let $ \{\phi_k(x), k = 1,2, \ldots \} $ be an orthonormal uniform 
boumded: $ \sup_{k,x}|\phi_k(x)| < \infty $ sequence of functions 
on some non - trivial measurable space $ (X, A, \mu) $ and  (in the $ L_2(X,\mu) $ sense) }
$$
f(x) = \sum_{k=1}^{\infty} c(k) \ \varphi_k(x).
$$
{\bf A).} {\it \ If $ c \in g(\psi, \nu), $ then }
$$
||f||L \left(N_1[\psi], X, \mu \right) \le C_9 \cdot( \max(1, \sup_{k,x} |\phi_k(x)|)) 
\cdot  ||c||g(\psi,\nu), \eqno(2.8)
$$
{\bf B.)} {\it Let $ c \in g(\alpha)  $ for some $ \alpha \in (0,1]. $ We assert that }
$$
||f|| L \left(N_{1/\alpha}, X, \mu \right) \le C_{10}(\alpha) \cdot (\max(1, \sup_{k,x} |\phi_k(x)|)
\cdot  \ ||c(\cdot)||g(\alpha). 
$$
{\bf Theorem 6.} {\it If $ f \in G(\alpha; \psi,\nu), \ $ where $ \alpha \ge 2, $  then }
$$
  \sup_{M \ge 1} ||S_M[f]||L \left(N^{(\alpha)}_d [\psi] \right) \le C_{11}(\alpha,\psi, N, \nu) \ 
||f||G(\alpha; \psi,\nu). \eqno(2.9)
$$

\section{ Auxiliary results}

{\bf Theorem 7.} {\it Let  $ N(u) = N(W,u) = \exp(W(\log u)), \ u > e^2,  \ 
\psi(p) = \exp(W^*(p)/p), \ p \ge 2, $ and $ X = T^d. $ 
We propose that the Orlicz's norm $ ||\cdot|| L(N) $ and the  norm 
$ ||\cdot||G(\psi) $ are equivalent.  Moreover, in this case  $ f \ne 0, \ f \in G(\psi) $ 
(or $ f \in L(N(W(\cdot),u) \ ) $ if and only if}  $ \ \exists C_{12},C_{13}, C_{14} \in 
(0,\infty) \ \Rightarrow \forall \ u > C_{14} $ 
$$
T(|f|,u) \le C_{12} \exp \left(-W \left(\log \left(u/C_{13} \right) \right) \right), \eqno(3.1)
$$
{\it where for each measurable function } $ f: X \to R $ 
$$
T(|f|,u) = mes \{x: |f(x)| > u \}.
$$

{\bf Proof } of theorem 7. \ A). Assume at first  that $ f \in L(N), \ f \ne 0.$ 
Without loss of generality we suppose that $ ||f||L(N) = 1/2. $ Then  
$$
\int_X N(W, |f(x)|) \ dx \le 1 < \infty.
$$
 The proposition (3.1) follows from Chebyshev's inequality such that in (3.1)
$ C_{12} = 1, \ C_{13} = C_{14} = 1/||f||L(N), \ f \ne 0. $ \par
B). Inversely, assume that $ f, \ f \ne 0 $ is a measurable function, $ f: X \to R^1 $ such 
that
$$
T(|f|,u) \le \exp(-W(\log u)), \ u \ge e^2.
$$
 We have by virtue of properties of the function $ W: $
$$
\int_X N(|f(x)|/e^2) \ dx = \int_{ \{x: |f(x) \le e^2\} } + \int_{ \{x: |f(x)| > e^2} =
I_1 + I_2;
$$
$$
I_1 \le \int_X N(1) \ dx = N(1),
$$
$$
I_2 \le \sum_{k=2}^{\infty} \int_{e^k < |f| \le e^{k+1} } 
\exp(W(|f(x)|/e^2)) \ dx \le
$$
$$
 \sum_{k=2}^{\infty} \exp((W(k-1)) \ T(|f|,k) \le
  \sum_{k=2}^{\infty} \exp(W(k-1) - W(k)) < \infty.
$$
 Thus, $ f \in L(N(W)) $ and 
$$
\int_X N(|f(x)|/e^2) \ dx \le N(1) + \sum_{k=2}^{\infty} \exp(W(k-1) - W(k)) < \infty.
$$
C). Let now $ f \in G(\psi); $ without loss of generality we can assume that $ ||f||G(\psi) 
=1. $ We deduce for $ p \ge 2:$ 
$$
\int_X |f(x)|^p \ dx \le \psi^p(p).
$$
We obtain using again the Chebyshev's  inequality:
$$
T(|f|,u) \le u^{-p} \psi^p(p) = \exp \left[-p \log u + p \log \psi(p) \right],
$$
and  after the minimization over $  p: \ u \ge \exp(2) \ \Rightarrow  $
$$
T(|f|,u) \le \exp \left( - \sup_{p \ge 2} (p \log u - p \log \psi(p)) \right) =
$$
$$
\exp \left( (p \log \psi(p))^*(\log u)  \right) = \exp(-W(\log u)).
$$
D). Suppose now that $ T(|f|,x) \le \exp(-W(\log x)), \ x \ge \exp(2). $ We conclude:
$$
\int_X |f(x)|^p \ dx = p \int_o^{\infty} x^{p-1} T(|f|,x) dx = p \int_0^{\exp(2)} 
x^{p-1} T(|f|,x) dx +
$$
$$
p \int_{\exp(2)}^{\infty} x^{p-1} T(|f|,x) dx \le p \int_0^{\exp(2)} x^{p-1} dx + p 
\int_{\exp(2)}^{\infty} x^{p-1} \ T(|f|,x) dx \le
$$

$$
 e^{2p} + \int_{\exp(2)}^{\infty} p x^{p-1} \exp(-W(\log x)) \ dx =
$$
$$
e^{2 p} + p \int_2^{\infty} \exp(py - W(y)) \ dy, \ p \ge 2.
$$
We obtain using Laplace's method and  theorem of Fenchel - Moraux ([8], p. 23 - 25):
$$
\int_X |f(x)|^p dx \le e^{2p} + C^p \exp \left(\sup_{y \ge 2} (py - W(y)) \right) = e^{2p} +
$$

$$
C^p \exp(W^*(p)) = e^{2p} + C^p \exp(p \log \psi(p)) \le C^p \ \psi^p(p).
$$
 Finally, $ ||f||G(\psi) < \infty. $ \par
For example, if $ m > 0, \ r \in R, $ then
$$
 f \in L \left(N_{m,r} \right) \Leftrightarrow  \sup_{p \ge 2} \left[|f|_p \ p^{-1/m} \
\log^{-r} p \right] < \infty \ \Leftrightarrow 
$$
$$
T(|f|,u)  \le C_0(m,r) \exp \left( - C(m,r) u^m \left( \log^{-mr} u \right) \right), \ u \ge 2.
$$

{\it Remark 3.} If conversely 
$$
T(|f|,x) \ge \exp(-W(\log x)), \ x \ge e^2,
$$
then for sufficiently large values of $ p; \ p \ge p_0 = p_0(W) \ge 2  $
$$
|f|_p \ge C_0(W) \ \psi(p), \ \  C_0(W) \in (0, \infty).\eqno(3.2)
$$

{\it Remark 4.} In this proof we used only the condition $ 0 < mes(X) < \infty. $ Therefore,
our  conclusions in theorem 7 are true in this more general case.\par

{\bf Theorem 8.} {\it Let $ \psi \in \Psi.$ We assert that $ f \in L^0(N[\psi]), $ or, 
equally, $ f \in G^0(\psi) $ if and only if }

$$
\lim_{p \to \infty} |f|_p /\psi(p) = 0. \eqno(3.3)
$$
{\bf Proof.} It is sufficient by virtue of theorem 7 to consider only the case of 
$ G(\psi) $ spaces.\par
 1. Denote $ G^{00}(\psi) = \{ f: \ \lim_{p \to \infty} |f|/\psi(p) = 0 \}. \ $ Let 
$ f \in G^0(\psi), \  f \ne 0. $ Then for arbitrary $ \delta = const > 0 $ there exists a 
constant $ B = B(\delta,f(\cdot)) \in (0, \infty) $ such that 

$$
||f-f I(|f| \le B) \ ||G(\psi) \le \delta/2.
$$
Since  $ |f| I(|f| \le B)| \le B, $ we deduce 
$$
|f I(|f| \le B)|_p/\psi(p) \le B/\psi(p).
$$
We obtain using triangular inequality  for sufficiently large values $ p: \ p \ge p_0(\delta)=
p_0(\delta,B) \ \Rightarrow $
$$
|f|_p /\psi(p) \le \delta/2 + B/\psi(p) \le \delta,
$$
as long as $ \psi(p) \to \infty $ at $ p \to \infty. $ Therefore $ G^0(\psi) 
\subset G^{00}(\psi). $\par
(The set $ G^{00}(\psi) $ is a closed subspace of $ G(\psi) $ with respect to the $ G(\psi) $ 
norm and contains all bounded functions.) \par
2. Inversely, assume that $ f \in G^{00}(\psi). $ We deduce denoting $ f_B = f_B(x) = 
f(x) I(|f|> B) $ for some $ B = const \in (0, \infty):$

$$
\forall Q \ge 2 \ \Rightarrow \lim_{B \to \infty} |f_B|_Q = 0. 
$$
 Further,
$$
||f_B||G(\psi) = \sup_{p \ge 2} |f_B|_p/\psi(p) \le \max_{p \le Q} |f_B|_p/\psi(p) +
$$

$$
\sup_{p > Q} |f_B|_p/\psi(p) \stackrel{def}{=} \sigma_1 + \sigma_2;
$$

$$
\sigma_2 = \sup_{p > Q} |f_B|_P/\psi(p) \le \sup_{p \ge Q} (|f|_p/\psi(p)) \le \delta/2
$$
for sufficiently large $ Q $ as long as $ f \in G^{00}(\psi). $  Let us now estimate the 
value $ \sigma_1: $
$$
\sigma_1 \le \max_{p \le Q} |f_B|_p/\psi(2) \le \delta/2
$$
for sufficiently large $ B = B(Q). $ Therefore,

$$
\lim_{B \to \infty} || f_B||G(\psi) = 0, \ \  f \in G^0(\psi).
$$
{\bf Theorem 9}. {\it Let $ \psi(\cdot) = \psi_N(\cdot), \theta(\cdot) = \theta_{\Phi}(\cdot) \ $ 
be a two functions on the classes $ \Psi $ with correspondent $ N \ - $ Orlicz's functions
$ N(\cdot), \Phi(\cdot): $

$$
N(u) = \exp \left \{\left[p \log \psi(p) \right]^*(\log u) \right\},
$$

$$
\Phi(u) = \exp \left \{ \left[ p \log \theta(p) \right]^*(\log u \right \}, \ \ u \ge \exp(2).
$$
We assert that  $ \lim_{p \to \infty} \psi(p)/\theta(p) = 0 $ if and only if $ N(\cdot) >> 
\Phi(\cdot).$ } \\
{\bf Proof } of theorem 9. A). Assume at first that $ \lim_{p \to \infty} \psi(p)/\theta(p) = 0. $
Denote $ \epsilon(p) = \psi(p)/\theta(p), $ then $ \epsilon(p) \to 0, \ p \to \infty. $ \par
 Let $ \{ f_{\zeta}, \ \zeta \in Z \} $ be arbitrary bounded in the $ G(\psi) $ sense set of 
functions:
$$
\sup_{\zeta \in Z} ||f_{\zeta}||G(\psi) = \sup_{\zeta \in Z} \sup_{p \ge 2} |f_{\zeta}|_p/\psi(p)
= C < \infty,
$$
then 
$$
\sup_{\zeta \in Z} |f_{\zeta}|_p/\theta(p) \le C \epsilon(p) \to 0, \ p \to \infty.
$$
 It follows from previous theorem that $ \forall \zeta \in Z \ f_{\zeta} \in G^0(\theta) $ and 
that the family $ \{f_{\zeta}, \ \zeta \in Z \} $ has uniform absolute continuous norm. Our 
assertion follows from lemma 13.3 in the book [6].\par
B). Inverse, let $ \Phi(\cdot) << N(\cdot). $ Let us introduce the measurable function 
$ f: X \to R $ such that $ \ \forall x \ge \exp(2) $

$$
\exp \left(-2 \left[p \log \psi(p) \right]^*(\log x) \right) \le T(|f|,x) \le
$$

$$
\exp \left( - \left[p \log \psi(p) \right]^*(\log x) \right).
$$
Then (see theorem 7) 
$$
f(\cdot) \in G(\psi), \ \ \ C_{15}(\psi) \ \psi(p) \le |f|_p \le C_{14}(\psi) \ \psi(p), \ \  p \ge 2.
$$
 Since $ f \in G(\psi), \  \Phi << N, $ we deduce that $ f \in G^0(\theta), $
and, following, 
$$
\lim_{p \to \infty} |f|_p/\theta(p) = 0.
$$
Therefore, $ \lim_{p \to \infty} \psi(p)/\theta(p) = 0. $ \par
{\bf Theorem 10}. {\it Let now $ X = R^d $ and $ \psi \in \Psi. $
We assert that the norms \\
$ ||\cdot|| L(N^{(\alpha)}, [\psi]) $ and $||\cdot|| G(\alpha, \psi), \ \alpha \ge 1 $ 
are equivalent.} \par
{\bf Proof.}
1. Let $ \forall p \ge \alpha \ \Rightarrow |f|_p \le \psi(p), \ f \ne 0. $  From 
Chebyshev inequality follows that
$$
\lim_{v \to \infty} T(|f|,v) = 0. 
$$
We can choose the values $ v, C_{17}, C_{18}, $ such that $ 0 < C_{17} < C_{18} < \infty, 
\ v \in (0,\infty) $
and
$$
 C_{17} \le T(|f|,v) \le C_{18} < \infty.
$$
Let us consider for some sufficiently small value $ \epsilon \in (0, \epsilon_0), \ 
\epsilon_0 \in (0,1) $ the following integral:
$$
I_{\alpha, N}(f) = \int_X N^{(\alpha)}(\epsilon |f(x)| \ dx = I_1 + I_2, 
$$
where
$$
I_1 = \int_{\{x: |f(x)| \le v \} } N^{(\alpha)}(\epsilon |f(x)|) \ dx,
\ I_2 = \int_{\{x: |f(x)|> v \} } N^{(\alpha)}(\epsilon |f(x)|) \ dx.
$$
Since for $ z \ge v $
$$
N^{(\alpha)}(z) \le C_{19}(\alpha, N(\cdot)) \ \cdot N( z),
$$
we have for the set $ X(v) = \{x, |f(x)| > v \}, $ using the result of theorem 7 for the 
space with finite measure:
$$
I_2 = \int_{X(v)} N^{(\alpha)}(\epsilon |f(x)| ) \ dx \le C_{20}(\alpha,N,\epsilon) \ 
||f||L(N^{(\alpha)}, X(v)) \le
$$

$$
C_{21} (\alpha, \epsilon, \psi) \sup_{p \ge \alpha}\left[ ||f||L_p(X(v))/\psi(p) \right] \le 
C_{21} \sup_{p \ge \alpha} |f|_p/\psi(p) < \infty.
$$
Further, since for $ z \in (0,v) \ \Rightarrow $
$$
N^{(\alpha)}(\epsilon z) \le C_{22}(v,\alpha, \epsilon) \ |z|^{\alpha},
$$
we have:
$$
I_1 \le C_{22}(\cdot) \int_X |f(x)|^{\alpha} \ dx < \infty.
$$
Thus, $ f \in L(N^{(\alpha)}[\psi]), \ ||f||L(N^{(\alpha)} [\psi]) < \infty.$ \par
2). We prove now the inverse inclusion. Let $ f \in L \left(N^{(\alpha)}[\psi] \right) $ and 
$$
||f||L \left(N^{(\alpha)}[\psi] \right) = 1.
$$
Hence for some $ \epsilon > 0 $
$$
\int_X N^{(\alpha)}(\epsilon |f(x)|) dx < \infty. 
$$
 It follows from the proof of
 theorem 7 and the consideration of two cases: $ |z| \le v; \ |z| > v $ the 
following elementary inequality: at $ p \ge \alpha $ and for all $ z > 0 \ \Rightarrow $
$$
|z|^p \le C_{23} (\alpha, \epsilon,N) \ N^{(\alpha)} (\epsilon |z|) \cdot \psi^p(p). 
$$
 We obtain for all values $ p, \ p \ge \alpha: $
$$
\int_{R^d} |f(x)|^p \ dx \le C_{24}^p(\alpha, \epsilon,\psi) \ \psi^p(p), \ \ \ ||f||G(\alpha;\psi) < 
\infty.
$$
 Note that theorems 7 - 10 are some generalizations of [1, p. 201 - 233], 
[15], p.29  and [8], section 1. The case of Hilbert transform of generalized functions $ \{ f \} $ 
is consider in [2].\\

\section{ Proofs of main results. }

 At first we consider the case Orlicz spaces, i.e. if the function $ f $ belongs to some 
exponential Orlicz space.\par
{\bf Proof of theorems 1,2.} Let $ X = [0, \ 2\pi]^d $ and $ f \in L(N[\psi]) $ for some 
$ \psi \in \Psi.$ Without los of generality we can assume that $ ||f||L(N[\psi]) = 1. $
From theorem 7 follows that 
$$
\forall p \ge 2 \ \Rightarrow \ |f|_p \le C_{25}(\psi) \ \psi(p).
$$
 From the classical theorem of M. Riesz (see [3], p.100 - 103; [9]) follows the inequality:
$$
\sup_{M \ge 1} |s_M[f]|_p \le K_1^d \ p^d \ \psi(p), \ \ K_1 = 2 \pi.
$$
It follows again from theorem 7  that 
$$
\sup_{M \ge 1} ||s_M[f] - f||L(N_d[\psi]) \le K_1^d + 1 < \infty.
$$
 For example, if $ N(u) = N_m(u) = \exp(|u|^m)-1 $ for some $ m = const \ge 1, $ then
$$
\sup_{M \ge 1} ||s_M[f] -f ||L(N_{m/(dm+1)}) \le C_{26}(d,m) \ ||f||L(N_m).  \eqno(4.1)
$$
 The "continual" analog of M.Riesz's inequality, namely, the case $ X = R, \ L(N) = 
L_p(R), \ p \ge 2: $
$$
\sup_{M \ge 2} |S_M[f]|_p \le K_2^d \ p^d \ |f|_p, \ \ K_2 = 1 
$$
 is proved, for example, in [10], p.187 - 188. [9]. 
This fact permit  us to prove also theorem 2. \par
{\bf Lemma 1.} {\it We assert that the "constant" $ m/(dm+1) $ in the estimation (3.5)
is exact. In detail, for all $ m \ge 1 $ there exists $ g = g_m(\cdot) \in L(N_m) $ such that 
$ \forall  \Delta \in (0,m) $ } 
$$
\sup_{M \ge 1} ||s_M[g]||L \left(N_{(m-\Delta)/(dm+1)} \right) = \infty.
$$
{\bf Proof} of lemma 1. It is enought to prove that
$$
\exists g \in L(N_m), \ \ ||H[g]||L(N_{(m-\Delta)/(dm+1)}) = \infty,
$$
where $ H[g] $ denotes the Hilbert transform on the $ [0, \ 2 \pi]^d, $
see [12], p. 193 - 197; [13]. Also it is enought to consider the case $ d = 1. $\par
 Let us introduce the function 
$$
g(x) = g_m(x) = | \log(x/(2 \pi))|^{1/m}. 
$$
Since for $ u > 0 $ 
$$
mes \{x: g_m(x) > u \} = \exp \left(- u^m \right),
$$
we conclude $ g_m(\cdot) \in L(N_m) \setminus L^0(N_m) $ (theorem 7).  
 Further, it is very simple to verify using the formula for  Hilbert transform ([3], p.192)  that
$$
C_{28}(m) \left(| \log (x/ (2 \pi))|^{(m+1)/m} + 1 \right) \le |H[g_m](x)| \le 
$$
$$
C_{29}(m) \left( |\log (x/(2 \pi))|^{(m+1)/m} + 1 \right).
$$
 Hence  $ \ \forall u \ge 2 $ 
$$
\exp \left( -C_{29}(m) u^{m/(m+1)} \right) \le mes \{ x, |H[g_m](x) > u \} \le 
$$
$$
\exp \left(-C_{30}(m) \ u^{m/(m+1)} \right).
$$
 It follows again from theorem 7 that 
$$
H[g_m] \in L(N_{m/(m+1)}) \setminus L(N_{m/(m+1)}). 
$$
Hence $ \forall \Delta \in (0,m) \ \Rightarrow \ H[g_m] \notin L(N_{(m-\Delta)/(m+1)}). $ \par
{\bf Proof} of theorem 3. Let us consider the following function:
$$
z(x) = z_L(x) = \sum_{n=8}^{\infty} n^{-1} \ L(n) \ \sin(nx). \eqno(4.2)
$$
 It is known from the properties of slowly varying functions ([14], p. 98 - 101) that 
the series (4.2) converge a.e. and at $ x \in (0, \ 2 \pi] $
$$
C_0 L(1/x) \le z(x) \le C L(1/x). 
$$
Therefore, at $ u \in [\exp (2), \infty) $ 
$$
  L^{-1}(C u)  \le  T(|z|,u) \le   L^{-1}(C_0 u). 
$$
 It follows from theorem 7 and (3.2)  that 
$$
z(\cdot) \in L(N) \setminus L^0(N), \ N(u) = L^{-1}(u), \ u \ge \exp(2).
$$
 From theorem 8 follows that 
$$
0 < C_0 \le |z|_p/\psi(p) \le C < \infty, \ p \ge 2, \ \psi(p) = \exp(W^*(p)/p). \eqno(4.3)
$$
 Note as a consequence that the series (4.2) {\it does not converge in the $ L(N) $ norm,}
as long as the system of functions $ \{\sin(nx) \} $ is bounded and hence in the case
when the series (4.2) {\it converge}  in the $ L(N) $ norm $ \Rightarrow z(\cdot) \in L^0(N). $   \par
 Let us suppose now that for some EOF $ \Phi(\cdot) $ with correspondence function $ \theta(p) $ 
the series (4.2) convergence in the $ L(\Phi) $ norm. Assume converse to the 
assertion of theorem 3, or equally that 
$$
\overline{\lim}_{p \to \infty} \theta(p)/\psi(p) > 0. \eqno(4.4)
$$
 Since the system of functions $ \{ \sin(nx) \} $ is bounded, $ z(\cdot) \in L^0(\Phi). $
By virtue of theorem 8 we conclude that
$$
\lim_{p \to \infty} |z|_p/\theta(p) = 0.
$$ 
 Thus, we obtain from (4.3)
$$
\lim_{p \to \infty} \psi(p)/\theta(p) = 0,
$$
in contradiction with (4.4). The cases $ X = [0, 2 \pi]^d, X = R^d $ are considered as well as 
the case $ X = [0, 2 \pi]. $ \par
 Now we consider the case when $ f \in G(a,b,\alpha, \beta).$ \par
{\bf Proof } of theorem 4. Let $ f \in G(1,b,\alpha,0), \ ||f||G(1,b,\alpha,0) = 1. $ 
Then $ \forall q \in (1,2] $
$$
|f|_q \le (q-1)^{-\alpha}. 
$$
We denote $ p = q/(q-1), $ then $ p \in [2, \infty). $ We will use the classical result of 
Hardy - Littlewood, Hausdorff - Young ([3], p.193; [16], p. 93;  [17], p. 26):
$$
|F[f]|_p \le (2 \pi)^{1/2} \ |f|_q.
$$
We have in our case:
$$
|F[f]|_p \le (2 \pi)^{1/2} \ (q-1)^{-a} = (2 \pi)^{1/2} (p/(p-1) - 1)^{-a} < (2 \pi)^{1/2} \ p^a.
$$
Our proposition it follows  from theorem 10.\\
{\bf Proof } of theorem 5. {\bf A.} We will use the classical result of  Paley and 
F.Riesz ([19], p.120).
Let $ \{ \phi_k(x), \ k = 1,2,\ldots \} $ be some orthonormal bounded sequence of functions.
 Then $ p \ge 2 \ \Rightarrow $
$$
|f|_p \le K_3 \cdot p \ \cdot \left( 1 + \sup_{k,x} |\phi_k(x)| \right) \cdot
 \left( \sum_k |c(k)|^p \left( |k|^{p-2} + 1 \right) \right)^{1/p}, \eqno(4.5)
$$
where $ K_3 $ is an absolute constant, $ f(x) = \sum_k c(k) \ \phi_k(x). $ \par
Let $ c(\cdot) \in g(\psi, \nu) $ and $ ||c||g(\psi,\nu) = 1. $ 
By definition of the $ g(\psi,\nu) $ norm 
$$
  \sum_{k = 1}^{\infty} |c(k)|^p \ \left(k^{p-2} + 1 \right)  \le \ ||c||^p g(\psi, \nu) \cdot  
\psi^p(p).
$$
Therefore
$$
|f|_p \le  K_3 \ C_{31} \ p \cdot \psi(p),
$$
and by virtue of theorem 7 $ \ f(\cdot) \in L(N_1[\psi]). $ \par
{\bf B.} Here we use the "discrete" inequality of Hausdorff - Young,  Hardy - Littlewood 
(see [16], p.101; \ [17], p.26:)
$$
|f|_p \le K_4 \ |c|_q, \ p \ge 2, \ q = p/(p-1), \ \ K_4 = 2 \pi.
$$
If $ ||c||g(\alpha) =1, $ then 
$$
|c|_q \le (q -1)^{\alpha},    |f|_p \le K_4 \ p^{\alpha}, \ p \ge 2. 
$$ 
Again from theorem 7 follows that $ f \in L \left(N_{1/\alpha} \right). $ \par
{\bf Proof } of theorem 6. The analog of inequality (4.5) in the case 
$$ 
F[f](t) = \int_R \exp(i t x) f(x) dx, \ d = 1,
$$
namely:
$$
|F[f]|_p \le K_5 \ p \ |f|_p(\nu),
$$
when $ f(\cdot) \in L_p(\nu) \subset G(\nu, \alpha, \psi), $ see, for example, in 
([16], p. 108).  Hence for all $ \ p \ge \alpha $
$$
||F[f]||L \left(N_d^{(\alpha)}[\psi] \right) \le K_5 \ ||f||G(\alpha; \psi, \nu). 
$$
(The generalization on the case $ d \ge 2 $ is evident). \par
 Note that the moment estimations for the wavelet transforms and Haar series are described
for example in the books [7], p. 21; [10], p.297  etc. It 
is easy to generalize our results on the cases Haar's or wavelet series and transforms.\par
 In detail, it is true in this cases the moment estimation for the partial sums (wavelet's or 
Haar's)
$$
|P_M[f]|_p \le K_6 \ |f|_p, \ X = [0,1], \ p \ge 1,
$$
where $ K_6 = 13 $ for Haar series on the interval $ X = [0,1] $  and $ K_6 = 1 $ for the classical 
wavelet series, in both the cases $ X = [0,1] $ and $ \ X = R.$
 Hence $ \forall \psi \in \Psi, \ f(\cdot) \in L(N[\psi]) $ 
$$
\sup_{M \ge 1} ||P_M[f] - f||L(N[\psi]) \le (K_6+1) \ ||f||L(N[\psi]). \eqno(4.6)
$$
 But (4.6) does not mean in general case the convergence 
$$
\lim_{M \to \infty} ||P_M[f] - f||L(N[\psi]) = 0, \eqno(4.7)
$$ 
as long as if (4.7) is true, then $ f(\cdot) \in L^0(N[\psi]) $ and conversely if 
$ f \in L^0(N[\psi]), $ then (4.7) holds.\par
 For the different generalizations of wavelet series the estimation (4.6) with constants 
$ K_6 $  not depending on $ p, \ p \ge 2 $ see, for example, in the books ([5], p. 221;
[18], p. 69). \\
{\bf Aknowledgement.} I am very grateful for support and attention to prof. V.Fonf (Beer - Sheva,
Israel).This investigation was partially supported by ISF (Israel Science Foundation) grant 
$ N^o $ 139/03. \\ 

\newpage

{\bf References} \\

1. Buldygin V.V., Mushtary D.I., Ostrovsky E.I, Pushalsky M.I. {\it New Trends in Probability 
Theory and Statistics.} Mokslas, 1992, Amsterdam, New York, Tokyo. \\
2. Chung Jaeyoung. {\it Hilbert Transform of Generalized Function of $L_p - $ Growth. }
Integral Transforms and Special Functions. 2001, $ N^o $ 2, p. 149 - 160.\\
3. Edwards R.E. {\it Fourier Series. A modern Introduction.} 1982, v.2; Springer Verlag, 1982. 
Berlin, Heidelberg, Hong Kong, New York. \\
4. Gord Sinnamon. {\it The Fourier Transform in Weighted Lorentz Spaces.} Publ. Math., 2003, 
v. 47, p. 3 - 29. \\
5. Hernandes E., Weiss G. {\it A First Course on Wavelets.} 1996, CRC Press, Boca Raton, New 
York.\\
6. Krasnoselsky M.A., Routisky Ya. B.{\it Convex Functions and Orlicz Spaces.} P. Noordhoff 
Ltd, 1961, Groningen. \\ 
7. Novikov I., Semenov E. {\it Haar Series and Linear Operators.} 1997, Kluwner Academic 
Publishers.  Boston, Dorderecht, London. \\
8. Ostrovsky E.I. {\it Exponential estimates for random fields and their 
 applications } (in Russian). 1999, OINPE, Obninsk.\\
9. Pichorides S.K. {\it On the best values of the constant in the theorem of 
M.Riesz, Zygmund and Kolmogorov.} Studia Math., 44, 1972, 165 - 179.\\
10. Pinsky M.A. {\it Introduction to Fourier Analysis and Wavelets.} Books Coll., 2002,
Australia, Canada.\\
11. Rao M.M., Ren Z.D. {\it Theory of Orlicz Spaces.} Marcel Dekker Inc., 1991. New York,
Basel, Hong Kong.\\
12. Rao M.M., Ren Z.D. {\it Applications of Orlicz Spaces.} Marcel Dekker Inc.,
2002. New York, Basel, Hong Kong.\\
13. Ryan R. {\it Conjugate Functions in Orlicz Spaces.} Pacific Journal Math., 1963, v. 13,
p. 1371 - 1377.  \\
14. Seneta E. {\it Regularly Varying Functions.} Springer Verlag, 1985.
Russian edition, Moscow,  Science, 1985.\\
15. Taylor M.E. {\it Partial Differential Equations. } v.3, Nonlinear Equations.
Springer Verlag, 1991. Berlin, Heidelberg, New York.\\
16. Titchmarsh E.C. {\it Introduction to the Theory of Fourier Integrals. }
Claredon Press, Second Edition. 1948. Oxford.\\
17. Wolff T.H. {\it Lectures on Harmonic Analysis.}
AMS, 2003, Providence, Rhore Island.\\
18. Wong M.W. {\it Wavelet Transforms and Localization Operators.} Birkhauser Verlag, 2002.
Berlin, Basel, Boston.\\
19. Zygmund A. {\it Trigonometrical Series.} Cambridge, University Press, 1968, V.2.\\
\end{document}